\documentclass[conference]{IEEEtran}
\IEEEoverridecommandlockouts

\usepackage{cite}
\usepackage{amsmath,amssymb,amsfonts}
\usepackage{algorithmic}
\usepackage{graphicx}
\usepackage{textcomp}
\usepackage{xcolor}
\def\BibTeX{{\rm B\kern-.05em{\sc i\kern-.025em b}\kern-.08em
    T\kern-.1667em\lower.7ex\hbox{E}\kern-.125emX}}
\begin{document}

\title{Optimal BESS Scheduling for Multi‐Market Participation in the Nordics \\
}

\author{\IEEEauthorblockN{ Zeenat Hameed, Chresten Træholt}
\IEEEauthorblockA{\textit{Technical University of Denmark, Department of Wind and Energy Systems} \\
}
}

\maketitle

\begin{abstract}
Battery energy storage systems (BESSs) can provide fast frequency reserves and energy arbitrage in Nordic electricity markets, but their limited energy capacity requires accurate revenue forecasts and coordinated bidding across multiple submarkets. This paper introduces a unified framework that employs generalized additive models (GAMs) to generate one‐week‐ahead forecasts of spot and FCR–N revenues in Denmark, Finland, and Norway using three years of hourly price and volume data. Forecast outputs feed a stochastic mixed‐integer optimizer that co‐optimizes BESS participation in FCR–N, FCR–D, spot markets, subject to state‐of‐charge constraints, inverter losses, and differing pay‐as‐bid and pay‐as‐clear rules. Comparative analyses evaluate forecast accuracy and quantify the impact of forecast errors on BESS bid acceptance, market selection, and profitability under realistic seasonal price patterns. Results demonstrate that spot markets exhibit consistently higher predictability than frequency markets, and that forecast errors modestly affect bid acceptance but do not alter overall market participation. The proposed approach provides BESS operators and investors with a  tool to assess revenue uncertainty and optimize multi‐market strategies in Nordic power systems.
\end{abstract}


\section{Introduction}

The Nordic region has committed to carbon neutrality by 2050, driving a rapid expansion of renewable energy sources (RES) such as wind and solar. While this transformation reduces emissions, it also introduces variability and low inertia into the grid, increasing demand for fast frequency reserves (FCR) and real‐time balancing \cite{b1,b6}. Battery energy storage systems (BESSs) are well suited to provide these services due to their fast response and optimal control, but their limited energy capacity requires accurate revenue forecasts to optimize bidding strategies \cite{b1,b2}.

In the Nordics, FCR–N, FCR–D, and spot markets operate under distinct pricing and procurement rules: frequency markets use pay‐as‐bid auctions, whereas spot markets clear on pay‐as‐clear \cite{b6}. Moreover, each TSO—Svenska Kraftnät (SE), Statnett (NO), Fingrid (FI), and Energinet (DK)—applies its own prequalification and market‐design requirements \cite{b4,b5}. Consequently, a BESS owner with assets in multiple countries must evaluate both spot and FCR revenue potential under differing market mechanisms \cite{b1,b3}.

This paper addresses three questions: (1) Can a unified model forecast spot and FCR–N revenues across Denmark, Finland, and Norway? (2) How does predictability compare between spot and frequency markets in each country? (3) To what extent do BESS bid strategies depend on forecast accuracy across multiple markets? We employ generalized additive models (GAMs) on three years of hourly price and volume data to capture non‐linear intraday and intraweek patterns \cite{b1,b3}. The resulting one‐week‐ahead forecasts feed a stochastic mixed‐integer BESS optimization that co‐optimizes bids in FCR–N, FCR–D, S–DCH (BESS sells energy), and S–CH (BESS buys energy) while enforcing state‐of‐charge limits, inverter losses, and market‐specific rules \cite{b2,b3}. This framework enables BESS operators to quantify revenue uncertainty and plan multi‐market strategies under Nordic market conditions.

\section{Market Mechanisms}
The Nordic power system—a single ENTSO‐E synchronous area except for Denmark’s DK1 (linked to Continental Europe)—is divided into bidding zones: Denmark (2), Sweden (4), Finland (1), and Norway (5). Its electricity market comprises day‐ahead and intra‐day energy trading (via Nord Pool) and regulation and balancing markets (managed by Svenska Kraftnät, Statnett, Fingrid, and Energinet) to ensure real‐time balance \cite{b4,b5,b6}. Day‐ahead and intra‐day markets set financial transactions for the next 24 hours, while the balancing market settles actual delivery and consumption \cite{b6}. Regulation and balancing markets procure services—FCR–N, FCR–D, FFR, aFRR, and mFRR—to maintain system stability. FCR–N, FCR–D, and FFR are power‐based auctions with availability payments per MW; FCR–N also pays per MWh. BESS units participate via droop control, adjusting output based on frequency deviations \cite{b1,b6}.

Though broadly aligned, Nordic regulation‐market mechanisms differ by TSO. TSOs jointly set total reserve needs; in 2024, these were FFR 300MW, FCR–D 1450MW, FCR–N 600 MW, and aFRR 300MW \cite{b6}. Denmark and Sweden share a joint FCR–N market with unrestricted exchange and constant volumes; Fingrid operates day‐ahead (pay‐as‐clear) and year‐ahead (fixed‐price) auctions; Statnett enforces droop control and procures hourly and two‐day‐ahead \cite{b6}. FCR–D follows the same procurement structure as FCR–N in each country \cite{b6}. FFR, introduced in 2020, is managed nationally by each TSO. For aFRR, Finland, Norway, and Sweden use pay‐as‐bid auctions; mFRR bids from all countries are pooled into a common activation market that opens 14 days before operation and clears at the marginal price \cite{b6}.
\begin{figure}
    \includegraphics[width=0.5\textwidth]{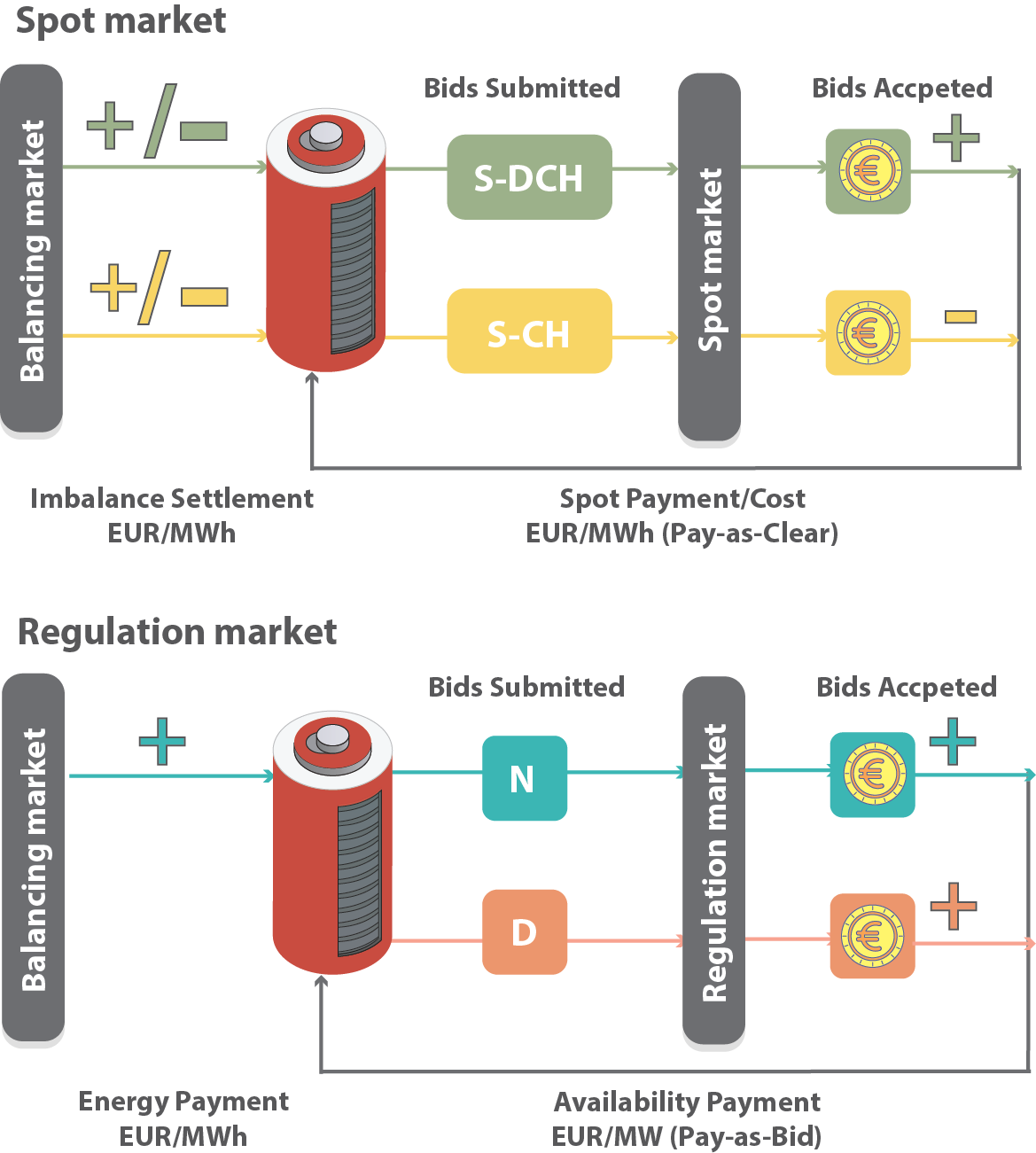}
    \caption{BESS payment mechanisms in FCR-N, FCR-D and spot markets}
    \label{fig:forecast_r2}
\end{figure}

\section{Method}

We forecast regulation‐market (RM) and energy‐market (EM) revenues using generalized additive models (GAMs), then optimize a battery energy storage system (BESS) across frequency and spot markets using those forecasts.

\subsection{GAM Revenue Forecasting}

GAMs model the expected response \(E(Y)\) as
\begin{equation}
g\bigl(E(Y)\bigr)
= \beta_{0}
  + \sum_{j=1}^{J} f_{j}(x_{j})\,,
\end{equation}
where each smoother
\begin{equation}
f_{j}(x_{j})
= \sum_{k=1}^{K} \beta_{j,k}\,b_{j,k}(x_{j})\,.
\end{equation}
The penalized log‐likelihood
\begin{equation}
\mathcal{L}_{p}
= \mathcal{L}
  - \lambda\,\boldsymbol{\beta}^{T} S\,\boldsymbol{\beta}
\end{equation}
balances fit \(\mathcal{L}\) against smoother complexity via penalty matrix \(S\) and smoothing parameter \(\lambda\).  Equation (1) defines the additive structure, (2) decomposes each smoother into basis functions \(b_{j,k}(\cdot)\), and (3) introduces a penalty \(\lambda\) to avoid overfitting.

For each country \(c\), define log‐revenues:
\begin{align}
FCRNrev_{h}^{c} &= \ln\bigl(FCRNprice_{h}^{c}\,FCRNvolume_{h}^{c}\bigr),\\
SPOTrev_{h}^{c} &= \ln\bigl(SPOTprice_{h}^{c}\,Consumption_{h}^{c}\bigr).
\end{align}
Equations (4)–(5) transform hourly prices and volumes into a log‐scale to stabilize variance.  Both GAMs then fit smooth functions in hour‐of‐day \((hours_{i})\), day‐of‐week \((days_{i})\), and their interaction:
\begin{equation}
\begin{split}
FCRNrev_{i}^{c}
&= f(hours_{i}) 
  + f(days_{i}) 
  + f(hours_{i}\times days_{i})\\
SPOTrev_{i}^{c}
&= f(hours_{i}) 
  + f(days_{i}) 
  + f(hours_{i}\times days_{i})
\end{split}
\end{equation}
Equation (6) shows the model form, where each smoother \(f(\cdot)\) captures non‐linear intraday and intraweek patterns; \(\varepsilon_{i}\) is Gaussian noise.  In R syntax:
\begin{equation}
\begin{split}
\texttt{gam}(Y) 
&= s(hours,\;bs=cr,\;k=24) \\
&\quad + s(days,\;bs=ps,\;k=7) \\
&\quad + ti(hours,\;days,\;k=c(24,7),\;bs=c(cr,ps)), \\
\text{family} &= \text{gaussian}.
\end{split}
\end{equation}
Equation (7) specifies cubic regression splines (\(bs=cr\)) for hours and penalized splines (\(bs=ps\)) for days, plus a tensor interaction (\(ti\)).  Forecast accuracy (with \(H=168\) hours per week) is measured by:

\begin{equation}
MAPE 
= \frac{1}{H}\sum_{h=1}^{H}\Bigl|\tfrac{A_{h}-F_{h}}{A_{h}}\Bigr|.
\end{equation}
Equation (8) defines mean absolute percentage error (MAPE).

\subsection{BESS‐Service Provision Optimization}

We define frequency and spot markets as:
\[
\begin{aligned}
M^{\mathrm{freq}} &= \{N, D\},\\
M^{\mathrm{spot}} &= \{\mathrm{S\text{-}DCH}, \mathrm{S\text{-}CH}\},\\
M &= M^{\mathrm{freq}} \cup M^{\mathrm{spot}}.
\end{aligned}
\]

Market‐power pairs are made to account for the MW bids that can be submitted in each market:
\[
MP = 
  \bigl(M^{\mathrm{freq}}\times TP^{\mathrm{freq}}\bigr)
  \cup 
  \bigl(M^{\mathrm{spot}}\times TP^{\mathrm{spot}}\bigr).
\]
Minute‐resolution frequency \(f_{st}\) and clearing prices \(P_{shm}^{\max}\) yield, via droop control,
\begin{equation}
E_{stmp}^{\mathrm{dch}},\quad
E_{stmp}^{\mathrm{ch}}
\quad
\forall (m,p)\in MP,\;s\in S,\;t\in T.
\end{equation}
Equation (9) defines market‐energy‐blocks (MEBs) of charge/discharge per minute for each scenario \(s\), period \(t\), market \(m\), and bid size \(p\).  The BESS state of charge \(z_{st}^{\mathrm{SOC}}\) satisfies:
\begin{align}
E^{\min} &\le z_{st}^{\mathrm{SOC}} \le E^{\max},\\
z_{s0}^{\mathrm{SOC}} &= 0.5\text{\,MWh},\\
z_{sT'}^{\mathrm{SOC}} &= 0.5\text{\,MWh},\\
ILF &= 0.10,
\end{align}
Equations (10)–(13) enforce minimum/maximum SOC limits \(E^{\min},E^{\max}\), set both initial and end‐of‐day SOC to 0.5 MWh, and define inverter‐loss factor \(ILF\).  Each hour \(h\in\{0,\ldots,23\}\) has \(N=60\) minutes and \(h_{t}=\lfloor t/N\rfloor\).

\subsubsection{Variables}
\[
\begin{aligned}
&x_{hmp}^{\mathrm{bid}}\in\{0,1\},\quad
x_{h}^{\mathrm{price}}\ge0,\quad
x_{shmp}^{\mathrm{bid,acc}}\in\{0,1\},\\
&z_{st}^{\mathrm{dch}},\,z_{st}^{\mathrm{ch}}\ge0,\quad
z_{st}^{\pm} = z_{st}^{\mathrm{dch}} - z_{st}^{\mathrm{ch}},\\
&z_{st}^{\mathrm{SOC}}\ge0,\quad
s_{st}^{\mathrm{dch}},\,s_{st}^{\mathrm{ch}}\ge0,\quad
w_{sh}^{\mathrm{ok}}\in\{0,1\},\\
&w_{shmp}^{\mathrm{avail}},\;
w_{sh}^{\mathrm{spot,dch}},\;
w_{sh}^{\mathrm{spot,ch}},\;
w_{sh}^{\mathrm{energy}}\ge0.
\end{aligned}
\]
Here, \(x_{hmp}^{\mathrm{bid}}\) indicates submitting a bid for hour \(h\), market‐power pair \((m,p)\).  \(x_{shmp}^{\mathrm{bid,acc}}\) equals 1 if that bid is accepted in scenario \(s\).  \(x_{h}^{\mathrm{price}}\) is the bid price (per MW).  \(z_{st}^{\mathrm{dch}},z_{st}^{\mathrm{ch}}\) are actual discharged/charged MWh in minute \(t\).  \(z_{st}^{\mathrm{SOC}}\) tracks SOC.  Slack variables \(s_{st}^{\mathrm{dch}},s_{st}^{\mathrm{ch}}\) record unmet discharge/charge.  \(w_{sh}^{\mathrm{ok}}=1\) if the BESS fulfills its energy commitment that hour.  Payments \(w_{shmp}^{\mathrm{avail}},w_{sh}^{\mathrm{spot,dch}},w_{sh}^{\mathrm{spot,ch}},w_{sh}^{\mathrm{energy}}\) follow from these commitments.

\subsubsection{Objective}
\begin{align}
\max \; &\sum_{s\in S} p_{s} \Bigl[
   \sum_{h,(m,p)\in MP^{\mathrm{freq}}} w_{shmp}^{\mathrm{avail}} \label{eq:obj}\\
&\quad + \sum_{h} w_{sh}^{\mathrm{spot,dch}} 
   - \sum_{h} w_{sh}^{\mathrm{spot,ch}} 
   + \sum_{h} w_{sh}^{\mathrm{energy}}
\Bigr].\notag
\end{align}
Equation (14) maximizes expected revenue: availability payments \(w_{shmp}^{\mathrm{avail}}\) from FCR‐N and FCR‐D, minus charge payments \(w_{sh}^{\mathrm{spot,ch}}\), plus discharge payments \(w_{sh}^{\mathrm{spot,dch}}\) and FCR‐N energy payments \(w_{sh}^{\mathrm{energy}}\).

\subsubsection{Constraints}

\paragraph{1) Bids:}
\begin{align}
\sum_{(m,p)\in MP} x_{hmp}^{\mathrm{bid}}
  &\le 1, 
  \label{eq:bid1}\\
  &\quad \forall\,h, \notag\\
x_{shmp}^{\mathrm{bid,acc}}
  &\le x_{hmp}^{\mathrm{bid}}, 
  \label{eq:bid2}\\
  &\quad \forall\,h, \notag\\
P_{shm}^{\max}\,x_{hmp}^{\mathrm{bid}} 
  &- \bigl(x_{h}^{\mathrm{price}} + \epsilon\bigr)
  \;\le\; M \, x_{shmp}^{\mathrm{bid,acc}}, 
  \label{eq:bid3}\\
  &\quad \forall\,(m,p)\in MP^{\mathrm{below}}, \notag\\
x_{h}^{\mathrm{price}} 
  &- \epsilon 
  - P_{shm}^{\max} 
  \;\le\; M \, x_{shmp}^{\mathrm{bid,acc}} 
  + M \,\bigl(1 - x_{hmp}^{\mathrm{bid}}\bigr), 
  \label{eq:bid4}\\
  &\quad \forall\,(m,p)\in MP^{\mathrm{above}}, \notag\\
bid_{\min} 
  &\le x_{h}^{\mathrm{price}} \le bid_{\max}, 
  \label{eq:bid5}\\
  &\quad \forall\,h. \notag
\end{align}

Equations \eqref{eq:bid1}–\eqref{eq:bid5} ensure:
\begin{itemize}
  \item \eqref{eq:bid1}: at most one bid per hour,
  \item \eqref{eq:bid2}: only placed bids can be accepted,
  \item \eqref{eq:bid3}: for “below” markets (FCR-N, FCR-D, S-DCH), bid price must not exceed the clearing price,
  \item \eqref{eq:bid4}: for “above” markets (S-CH), bid price must be at least the clearing price,
  \item \eqref{eq:bid5}: bid prices lie within \([bid_{\min},\,bid_{\max}]\).
\end{itemize}

\paragraph{2) Energy Flows \& SOC:}
\begin{align}
z_{stm}^{\mathrm{dch}} + s_{stm}^{\mathrm{dch}}
  &= E_{stm}^{\mathrm{dch}}\,x_{sh_{t}m}^{\mathrm{bid,acc}}, 
  \nonumber\\
  &\quad \forall\,s,t,(m,p),\\
z_{stm}^{\mathrm{ch}} + s_{stm}^{\mathrm{ch}}
  &= E_{stm}^{\mathrm{ch}}\,x_{sh_{t}m}^{\mathrm{bid,acc}}, 
  \nonumber\\
  &\quad \forall\,s,t,(m,p),\\
z_{st}^{\pm}
  &= \frac{1}{1-ILF}\,z_{st}^{\mathrm{dch}}
    - (1-ILF)\,z_{st}^{\mathrm{ch}}, 
  \nonumber\\
  &\quad \forall\,s,t,\\
E^{\min} \;\le\; z_{st}^{\mathrm{SOC}} - z_{st}^{\pm}
  &\;\le\; E^{\max}, 
  \nonumber\\
  &\quad \forall\,s,t,\\
z_{s0}^{\mathrm{SOC}} &= 0.5,\quad
 z_{sT'}^{\mathrm{SOC}} = 0.5, 
  \nonumber\\
  &\quad \forall\,s,\\
z_{st}^{\mathrm{SOC}}
  &= z_{s,t-1}^{\mathrm{SOC}} - z_{s,t-1}^{\pm}, 
  \nonumber\\
  &\quad \forall\,s,t\ge1,\\
\sum_{t:h_{t}=h} \bigl(s_{st}^{\mathrm{dch}} + s_{st}^{\mathrm{ch}}\bigr)
  &\le M \,(1 - w_{sh}^{\mathrm{ok}}), 
  \nonumber\\
  &\quad \forall\,s,h.
\end{align}

Equations (20)–(26) set per‐minute energy flows and SOC dynamics: (20)–(21) tie actual and slack charge/discharge to the MEBs \(E_{stm}^{\mathrm{dch}},E_{stm}^{\mathrm{ch}}\) and accepted bids; (22) computes net signed energy \(z_{st}^{\pm}\) after accounting for inverter‐loss factor \(ILF\); (23) enforces SOC limits; (24) fixes initial/final SOC; (25) updates SOC over time; (26) ensures that if \(w_{sh}^{\mathrm{ok}}=1\) (fulfilled), no slack is used, otherwise slack can exceed zero.

\paragraph{3) Market Payments:}

\emph{FCR‐N and FCR-D Availability:}
\begin{align}
w_{shmp}^{\mathrm{avail}}
  &\le p\,x_{h}^{\mathrm{price}},\\
w_{shmp}^{\mathrm{avail}}
  &\le M\,x_{shmp}^{\mathrm{bid,acc}},\\
w_{shmp}^{\mathrm{avail}}
  &\ge p\,x_{h}^{\mathrm{price}}
       - M \,(1 - x_{shmp}^{\mathrm{bid,acc}}).
\end{align}
Equations (27)–(29) linearize the product \(p\,x_{h}^{\mathrm{price}}\,x_{shmp}^{\mathrm{bid,acc}}\) using big‐\(M\), so that if \(x_{shmp}^{\mathrm{bid,acc}}=1\), the availability payment equals \(p\,x_{h}^{\mathrm{price}}\); otherwise it is zero.

\emph{Spot‐Discharge (S‐DCH):}
\begin{equation}
w_{sh}^{\mathrm{spot,dch}}
= \sum_{\substack{(m,p)\in MP^{\mathrm{spot}}\\ m=\mathrm{S\text{-}DCH}}}
  p \, P_{shm}^{\max} \, x_{shmp}^{\mathrm{bid,acc}}.
\end{equation}
Equation (33) computes pay‐as‐clear discharge revenue: each accepted bid \((m=\text{S‐DCH},p)\) yields \(p\cdot P_{shm}^{\max}\).

\emph{Spot‐Charge (S‐CH):}
\begin{equation}
w_{sh}^{\mathrm{spot,ch}}
= \sum_{\substack{(m,p)\in MP^{\mathrm{spot}}\\ m=\mathrm{S\text{-}CH}}}
  p \, P_{shm}^{\max} \, x_{shmp}^{\mathrm{bid,acc}}.
\end{equation}
Equation (34) computes cost of charging energy on the spot “charge” submarket (also pay‐as‐clear).

\emph{FCR‐N Energy:}
\begin{align}
w_{sh}^{\mathrm{energy}}
  &\;\le\;
  \sum_{t:h_{t}=h} \bigl[
    C_{sh}^{\mathrm{up}}\,z_{stm}^{\mathrm{dch}}
    + C_{sh}^{\mathrm{down}}\,z_{stm}^{\mathrm{ch}}
  \bigr],\\
w_{sh}^{\mathrm{energy}}
  &\;\le\; M \, w_{sh}^{\mathrm{ok}},\\
w_{sh}^{\mathrm{energy}}
  &\;\ge\;
  \sum_{t:h_{t}=h} \bigl[
    C_{sh}^{\mathrm{up}}\,z_{stm}^{\mathrm{dch}}
    + C_{sh}^{\mathrm{down}}\,z_{stm}^{\mathrm{ch}}
  \bigr]
  - M\,\bigl(1 - w_{sh}^{\mathrm{ok}}\bigr).
\end{align}
Equations (35)–(37) linearize the product of per‐minute FCR‐N energy prices \(C_{sh}^{\mathrm{up}},C_{sh}^{\mathrm{down}}\) with actual discharged/charged energy \(z_{stm}\), gated by \(w_{sh}^{\mathrm{ok}}=1\) to ensure energy payments only accrue if hourly net balance is nonnegative.

\paragraph{4) Variable Domains:}
\[
\begin{aligned}
&x_{hmp}^{\mathrm{bid}}, \; x_{shmp}^{\mathrm{bid,acc}}, \; w_{sh}^{\mathrm{ok}}
  \;\in\; \{0,1\},\\
&x_{h}^{\mathrm{price}}, \; z_{st}^{\mathrm{SOC}}, \; z_{st}^{\mathrm{dch}}, \; z_{st}^{\mathrm{ch}}, \; s_{st}^{\mathrm{dch}}, \; s_{st}^{\mathrm{ch}},\\
&w_{shmp}^{\mathrm{avail}}, \; w_{shmp}^{\mathrm{avail,D}}, \; w_{sh}^{\mathrm{spot,dch}}, \; w_{sh}^{\mathrm{spot,ch}}, \; w_{sh}^{\mathrm{energy}}
  \;\ge\; 0.
\end{aligned}
\]
Equation (38) enumerates binary and nonnegative domains for all decision and auxiliary variables.

\section{Results}

We trained six GAMs—two per country (Denmark, Finland, Norway)—to forecast weekly revenues in the spot and FCR‐N markets using rolling two‐week windows of hourly data. In‐sample fit (adjusted \(R^2\)) and out‐of‐sample accuracy (MAPE) were:

\medskip
\noindent\textbf{Spot market:}
\[
\begin{aligned}
R^2_{\mathrm{adj}} &= \{0.77\ (\mathrm{DK}),\,0.67\ (\mathrm{FI}),\,0.79\ (\mathrm{NO})\},\\
\mathrm{MAPE} &= \{0.91\%,\,0.74\%,\,0.94\%\}.
\end{aligned}
\]

\noindent\textbf{FCR–N:}
\[
\begin{aligned}
R^2_{\mathrm{adj}} &= \{0.92\ (\mathrm{DK}),\,0.46\ (\mathrm{FI}),\,0.06\ (\mathrm{NO})\},\\
\mathrm{MAPE} &= \{1.52\%,\,2.26\%,\,4.34\%\}.
\end{aligned}
\]
Spot‐market smoothers (hourly, daily, interaction) were nearly identical across countries; FCR–N smoothers showed a midday/midweek dip only in Denmark.

\begin{figure}[h]
    \centering
    \includegraphics[width=0.5\textwidth]{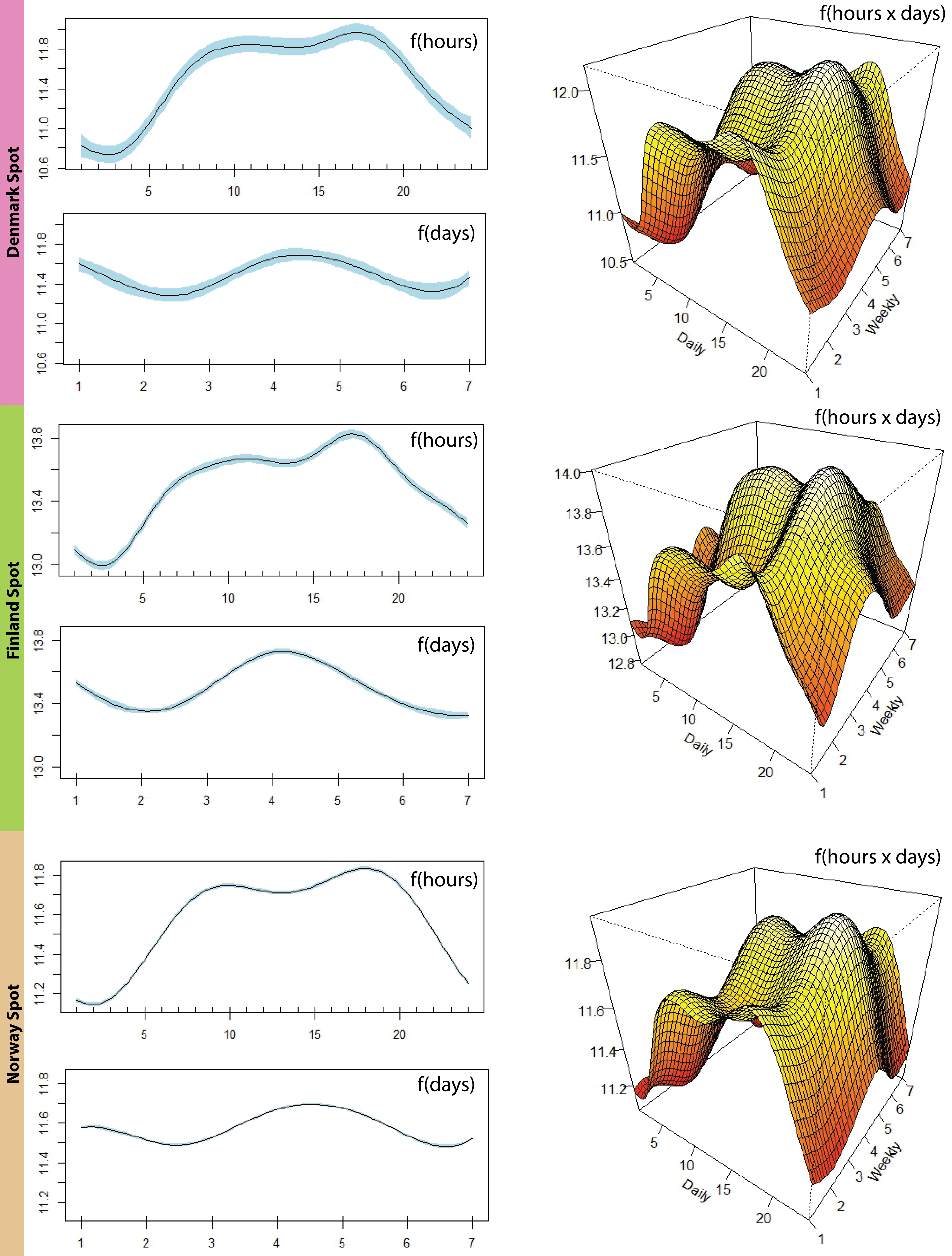}
    \caption{Spot‐market smoothers (hourly, daily, interaction) by country.}
    \label{fig:smoothers_spot}
\end{figure}

\begin{figure}[h]
    \centering
    \includegraphics[width=0.5\textwidth]{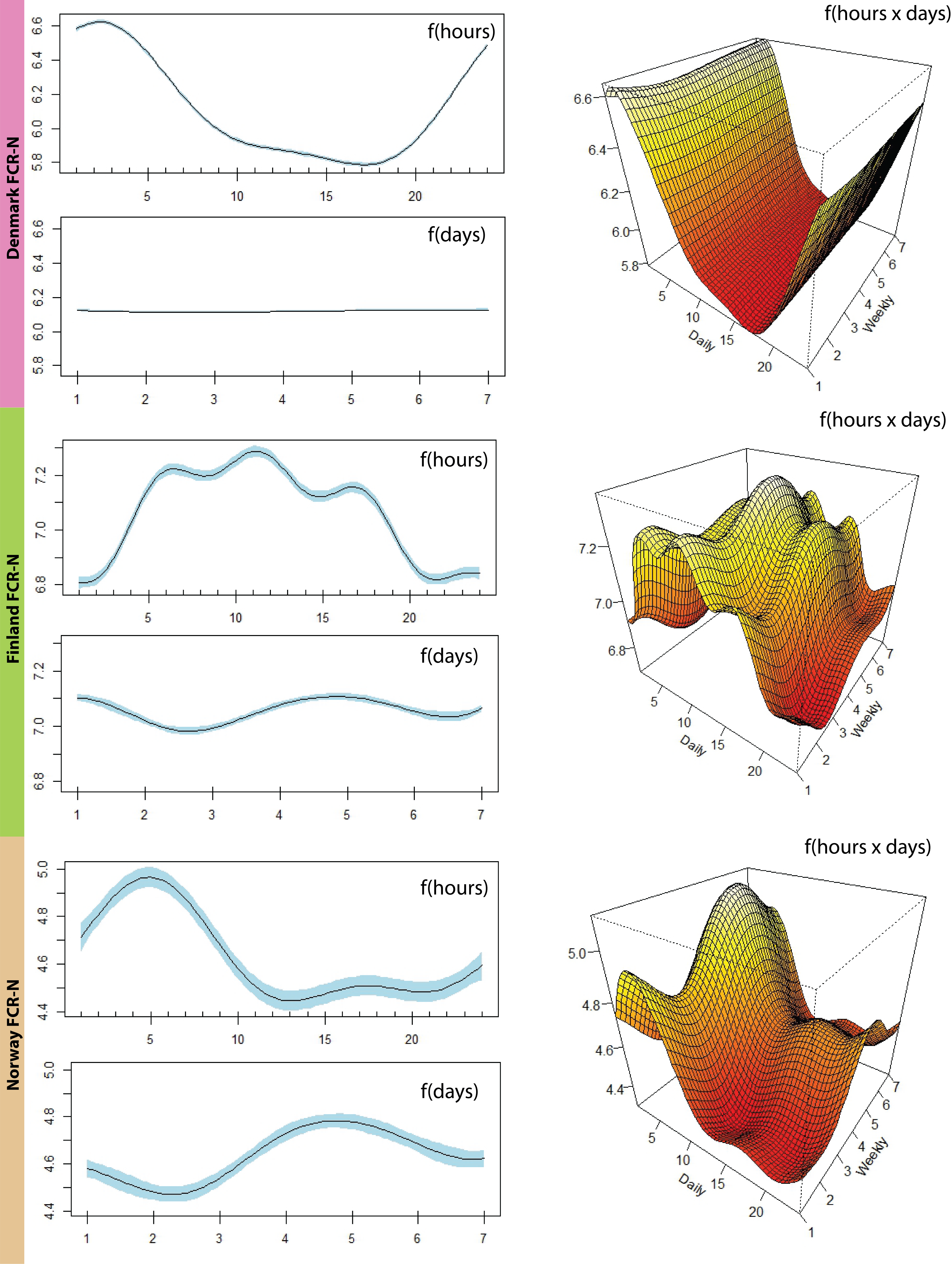}
    \caption{FCR–N smoothers by country.}
    \label{fig:smoothers_fcrn}
\end{figure}

Over 2019–2021, weekly MAPE distributions in three ranges (0–5\%, 5–10\%, 10–15\%) showed:
\begin{itemize}
    \item Denmark: 60.65\% of spot weeks and 87.1\% of FCR–N weeks had MAPE \(\le5\%\).
    \item Finland \& Norway: 85.16\% of spot weeks had MAPE \(\le5\%\), but only 53.55\% (FI) and 7.1\% (NO) of FCR–N weeks did.
\end{itemize}

\begin{figure}[h]
    \centering
    \includegraphics[width=0.5\textwidth]{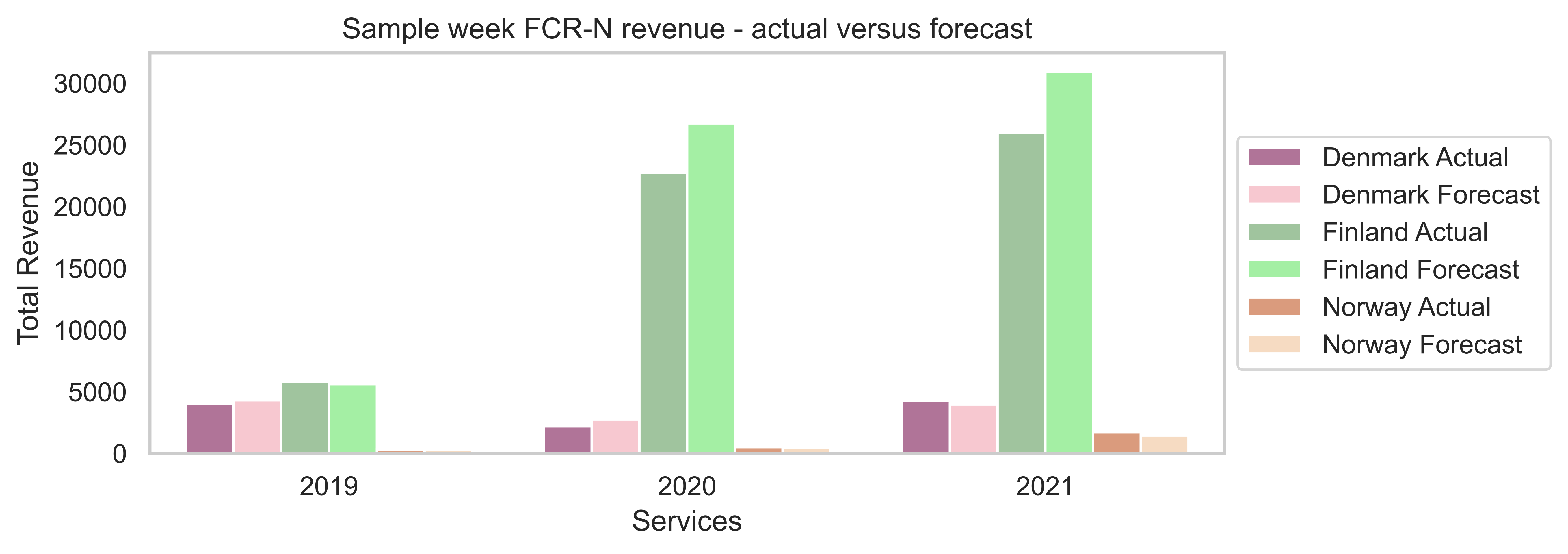}
    \caption{Actual versus forecasted valued for FCR–N, 2019–2021.}
    \label{fig:mape_distribution}
\end{figure}

Figures~\ref{fig:smoothers_spot} and~\ref{fig:smoothers_fcrn} illustrate the estimated smoothers for the spot and FCR–N markets across Denmark, Finland, and Norway. Specifically, Fig.~\ref{fig:smoothers_spot} shows the hourly, daily, and interaction smoothers for the spot market in each country, while Fig.~\ref{fig:smoothers_fcrn} displays the corresponding FCR–N smoothers. Figure~\ref{fig:mape_distribution} compares actual versus forecasted FCR–N values for each country.

\begin{figure}[h]
    \centering
    \includegraphics[width=0.5\textwidth]{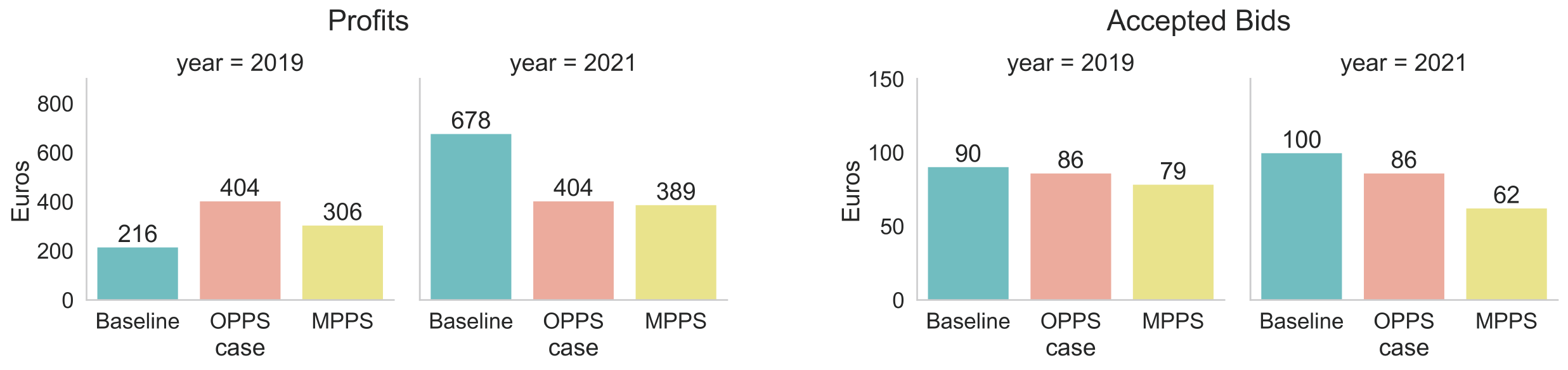}
    \caption{Estimated BESS profits and bid acceptance hours using original, GAM, and AI forecasts for 2019 and 2021.}
    \label{fig:sensitivity_volatility_profits_bids}
\end{figure}

To evaluate forecast impact on profits, we compared original data with GAM‐based forecasts (\(W_{\mathrm{rep}}^{2019},W_{\mathrm{rep}}^{2021}\)) and AI‐based forecasts (\(S_{\mathrm{rep}}^{2019},S_{\mathrm{rep}}^{2021}\)). Figure~\ref{fig:sensitivity_volatility_profits_bids} shows BESS profits (left axis) and bid acceptance hours (right axis) for original (teal), GAM (orange) illustrated as one-price-point-scenario (OPPS), and AI (yellow) illustrated as multiple-price-points-scenario (MPPS):

\begin{itemize}
    \item \textbf{2019:} GAM‐based profits were 30\% higher and AI‐based 17\% higher than original, while bid acceptance hours fell by 4\% (GAM) and 11\% (AI).
    \item \textbf{2021:} GAM‐ and AI‐based profits were 25\% and 27\% lower than original, with bid acceptance down by 4\% (GAM) and 38\% (AI).
\end{itemize}
In all scenarios, the optimizer predominantly chose \(m = N\). Higher optimal bid prices under forecasted data stemmed from misestimated pay‐as‐bid availability payments, inflating revenue despite fewer accepted bids.

We then compared BESS performance under a multi‐bid scenario—three discrete power levels \(\{0.9,0.6,0.3\}\,\mathrm{MW}\) for FCR–N (markets \(N,D\)) and three energy levels \(\{0.8,0.6,0.4\}\,\mathrm{MWh}\) for spot (S–DCH, S–CH)—against a single‐bid, flexible‐SOC baseline. In a representative 2019 week, multi‐bid profits were 2\% lower and costs 3.5\% lower; in 2021, profits were 0.8\% lower. Thus, multi‐bidding has negligible impact on net outcome.

\begin{figure}[h]
    \centering
    \includegraphics[width=0.5\textwidth]{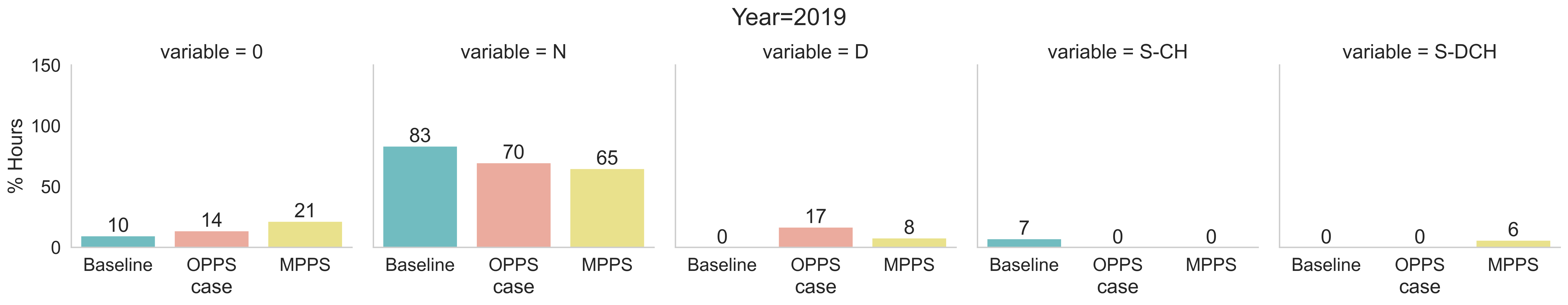}
    \caption{Profit and cost differences between multi‐bid and single‐bid strategies for representative weeks.}
    \label{fig:sensitivity_volatility_markets}
\end{figure}

Under a fixed‐SOC constraint (start/end each day at 0.5\,MWh), the BESS shifted capacity into FCR–N: 2019 participation rose from 71\% to 83\% of hours; 2021 from 4\% to 17\%. FCR–N’s minute‐by‐minute requirements accommodate fixed SOC, while spot’s block commitments risk SOC violations if SOC cannot vary. Finally, seasonal price patterns shaped market allocation and profitability. In 2019, midday and midweek FCR–N spikes led the BESS to spend 63\%–83\% of hours in \(N\), 0\%–15\% in spot discharge, 4\%–7\% in spot charge, and 0\% idle. In 2021, FCR–D spikes caused 83\%–96\% of hours in \(D\), no spot participation, and idle hours fell from 10\%–19\% to 0\%–1\%. Consequently, the higher 2021 price spreads generated approximately 49\% greater profit than the same week in 2019.

\section{Conclusions}

In Finland and Norway, spot‐market revenues are more predictable (higher adjusted \(R^2\), lower MAPE) than FCR–N, due to consistent hourly, daily, and interaction patterns. In Denmark, FCR–N revenues are more predictable than in FI/NO and even exceed Danish spot‐market predictability. GAM forecast errors slightly reduce bid acceptance hours but do not change market‐choice patterns. Multi‐bid options have negligible impact on profit and cost, so single‐bid strategies suffice. A fixed daily SOC shifts BESS participation toward FCR–N, as its minute‐by‐minute requirements allow more frequent rebalancing than spot. Price seasonality drives behavior: in 2019, high midday/midweek FCR–N prices led to predominantly FCR–N hours; in 2021, elevated FCR–D prices caused nearly exclusive FCR–D participation and reduced idle hours from ~15\% to 1\%. As a result, 2021’s price dynamics produced ~49\% higher profit than the same week in 2019.

\end{document}